\documentclass[a4paper,11pt,twoside]{amsart}
\usepackage[latin1]{inputenc}
\usepackage[T1]{fontenc}
\usepackage[french]{babel}
\usepackage{amssymb,pifont}
\usepackage{pstcol,pstricks,pst-node,pst-tree}
\usepackage{graphicx}
\renewcommand{\bf}{\bfseries}
\renewcommand{\it}{\itshape}

\input xy
\xyoption{all}

\begin{document}
{\centering \large \bf Geometry of certain Lie-Frobenius groups}\\
\vskip 2mm
{\bf\it Jean Michel Dardi\'e, Alberto Medina and Hass\`ene Siby}

\vskip 1cm

\begin{minipage}[9pt]{12cm}
{\bf Abstract.}
Let be $G_{n,p}({\mathbb K})= M_{n,p}(\mathbb K) \rtimes GL({\mathbb K}^n)$, $\mathbb K=\mathbb R$ or $\mathbb C$, the semi-direct product of the additive group of matrices $M_{n,p}(\mathbb K)$ by the group $GL({\mathbb K}^n) $ where the action is done by multiplication of matrices. If $n=kp$,$p\geq 1$ the Lie group $G_{n,p}({\mathbb K})$ admits an exact left invariant symplectic form [{\bf 10}]. We study the geometry of this symplectic manifold. If $n>p$ ( resp. $n=p$ ) we prove that $G_{n-p,p}({\mathbb K})$ (resp $G_{n-1,1}({\mathbb K})$ ) is the symplectic reduction of the symplectic orthogonal $(M_{n,p}(\mathbb K))^{\bot}$ of $M_{n,p}(\mathbb K)$ in $G_{n,p}({\mathbb K})$ and reciprocally that $G_{n,p}({\mathbb K})$ is a symplectic double extension, in the sense of [{\bf 6}], of $G_{n-p,p}({\mathbb K})$ ( resp. of $G_{n-1,1}({\mathbb K})$). Moreover we show that $G_{n,p}({\mathbb K})$ admits two left transverse Lagrangian foliations (with affine and closed leaves) . Consequently there exists (a left invariant) a canonical torsion free symplectic connection on $G_{n,p}({\mathbb K})$.
\end{minipage}
\vskip 1cm
{\bf\large  Introduction} 
\vskip 1mm
The group $G_{n,1}({\mathbb K}) = Aff({\mathbb K}^n)$ where $\mathbb K=\mathbb R$ or $\mathbb C$ admits left invariant symplectic structures which are all exact because $H^2(aff({\mathbb K}^n), \mathbb K)=0 $([{\bf 9}]).\\
To give an invariant symplectic form on $G_{n,1}({\mathbb K})$ is to give a linear form $\alpha$ on the Lie algebra $aff({\mathbb K}^n)$ such that the 2-cobound $\delta\alpha$ is non degenerate. That is to say $\alpha$ is a point with an open orbit under the coadjoint action of $G_{n,1}({\mathbb K})$. These points have been characterized in [{\bf 9}] within the framework of study of affine group representations and in [{\bf 3}] where is carried on the study of left invariant symplectic structures of affine group initiated in [{\bf 4}]. In [{\bf 10}], it showed that coadjoint action of Lie group semi-direct product of $M_{n,p}(\mathbb K)$ and $GL({\mathbb K}^n) $ by matricial product, where $M_{n,p}(\mathbb K)$ indicate the additive group of $(n\times p)$-matrices, admits open orbits if only if $n=kp$. We extend to these groups denoted $G_{n,p}({\mathbb K})$ the result obtained in [{\bf 3}] for the classical affine group: existence of a unique symplectic structure up isomorphism (Theorem 2.7), existence  of transverse symplectic foliations (Theorem 2.6), as well as of a Lagrangian bi-foliation with closed leaves (Theorem 3.2). Such a pair of Lagrangian foliations is important in the quantization procedure (polarization) and implies the existence of a canonical (torsion free) symplectic connection on $G_{n,p}({\mathbb K})$. The natural left action of $G_{n,p}({\mathbb K})$ on itself being hamiltonian we can provide (Theorem 2.1) a similar of structure's theorem of simply connected symplectic Lie group from [{\bf 5}].\\
In the third part we study up the fibrations of this theorem for explicit suggest a construction of $ \left( G_{n,p}({\mathbb K}),d\alpha_{1}^{+}\right) $ when the generalized double extension of [{\bf 5}] don't apply in this case.\\
To complete this introduction recall that a symplectic Lie group $(G, \omega^{+})$ has an affine structure given by the left invariant connection ( torsion free and zero curvature ) $\nabla$ defined for all $x, y, z \in T_{\varepsilon}(G)$ by :
$$\omega^{+} (\nabla^{+}_{x^{+}}y^{+}, z^{+}) = - \omega^{+} (y^{+}, [x^{+}, z^{+}])$$
where $x^{+}$ denotes the left invariant vector field associated to $x$.\\
Then one says that the pair $(G, \nabla^{+})$ is an affine Lie group. In this case the product on $\mathcal G$, $xy = (\nabla^{+}_{x^{+}}y^{+})$ is with associator left invariant and verifies the condition $xy - yx = [x, y]$. This structure plays a central rule in this study.
\vskip 2mm 
{\bf 1. Left invariant symplectic structures on the Lie groups $G_{n,p}$}
\vskip 0.2mm
In what follows $G_{n}({\mathbb K})$ denotes the classical affine group where $\mathbb K=\mathbb R$ or $\mathbb C$ and  ${\mathcal G}_{n}({\mathbb K}) $ is its Lie algebra. By analogy, we denote $ G_{n,p}({\mathbb K})$ the group semi-direct product $M_{n,p}(\mathbb K) \rtimes Gl({\mathbb K}^n)$ with $n=kp, p\geq 1$ and ${\mathcal G}_{n,p}({\mathbb K}^n)$, its Lie algebra. Obviously $ G_{n,1}({\mathbb K})$ is isomorphic to $Aff({\mathbb K}^n)$. Considering that $H^2(aff({\mathbb K}^n), \mathbb K)=0 $ (see [{\bf 3}]) any left invariant alternate 2-form is invariantly exact. This result becomes general to ${\mathcal G}_{n,p}({\mathbb K})$ in the following way:\\
{\bf 1.1 Lemma.} Every left invariant symplectic form on $G_{n,p}(\mathbb K)$ is exact.\\
\vskip 1mm
{\bf Proof.} If $k=p=1$; ${\mathcal G}_{1,1}$ is isomorphic to $aff({\mathbb K}) $  and the result follows.\\
Assume that $k$ (or $p$) is greater than 1.
 Let be $\omega \in Z^2({\mathcal G}_{n,p};\mathbb K) $, ${\mathcal G}_{n,p}=M_{p,n}(\mathbb K) \rtimes gl(n)$. This means that we have 
$$\oint \omega ([a,b],c)=0 \qquad (\mbox{\bf *})  \quad \, \mbox {for every  $a$,$b$,$c$ in \,}{\mathcal G}_{n,p}.$$\\
If we take $a=(x,u)$, $b=(y,v)$ and $c=(0,Id)$\\
({\bf *}) implies $\omega (x,y)=0$ for every $x,y \in M_{n,p}$, and consequently $\omega (x,b) + \omega(a,y) = \omega(I,uy-vx+[u,v])$.\\
But this is equivalent to 
$$\omega (a,b) = \omega(I,ux-vy)$$
$$ = -\delta\beta (a,b) \quad \mbox{ with \, } \beta = \omega (I,\cdot ) \mbox{  \qquad $\blacksquare$}$$

{\bf  N.B} In fact for all Lie algebra ${\mathcal G}$ having an element $a \in \mathcal G$ such that $ad_a$ is projector, we have $H^2({\mathcal G},\mathbb K)=0 $.\\
{\bf Remark 1.} The mapping
$$ M_{p,n}(\mathbb K) \times gl(n)  \longrightarrow (M_{n,p}(\mathbb K) \times gl(n))^{*}$$
$$(H,M) \longmapsto \alpha_{(H,M)}$$ 
given by 
$$ \alpha_{(H,M)}(N,V)=tr(N.H)+tr(M.V)$$
is a linear isomorphism.\\
In the following, the dual space ${\mathcal G}_{n,p}^{*}$ of ${\mathcal G}_{n,p}$ is identified with $M_{p,n}(\mathbb K) \times gl(n)$.
To give an invariant closed $2$-form on $G_{n,p}(\mathbb K)$ it is equivalent to give a linear form $\alpha$ on ${\mathcal G}_{n,p}({\mathbb K})$ . The $2$-form $d\alpha^+$, where $\alpha^+$ denotes the left invariant $1$-form defined by $\alpha$ is symplectic if and only if $\alpha$ has an open orbit under the coadjoint representation.
\vskip 1mm
Using the Remark 1 and the natural embedding of ${\mathcal G}_{n,p}$ in $GL({\mathbb K}^{n+p}) $ we can show
\vskip 2mm
{\bf 1.2 Lemma.} The coadjoint representations of ${\mathcal G}_{n,p}$ and $ G_{n,p}$ are given by the following formulas: 
$$ ad^*_{(0,u)}(H, N)=(-H, [u,N]) \leqno (i)$$
$$ad^*_{(x,0)}(H,N)=(0,x.H) \leqno(ii)$$
$$Ad^*_{(0,U)}(H,N)=(H.U^{-1}, UNU^{-1}) \leqno(iii)$$
$$Ad^*_{(X,Id_{{\mathbb K}^n})}(H,N)=(H, N+X.H) \leqno(iv)$$
This implies
$$ Ad^*_{(X,U)}(H,N)=(H.U^{-1}, UNU^{-1}+XHU^{-1}) $$
Consequently the coadjoint orbit $Ad^{*}_{G_{n,p}}(\alpha)$ of $\alpha =(H_0,N_0)$ with $H_{0}=(0, I_{p}) \in M_{p,n}$ and 
$$N_{0}=\left(\begin{array}{ccc} 
 0 & 0 & 0\\
 I_p & 0 & 0 \\
 0 & I_p & 0
\end{array}\right)$$ is open.
Denote by  $H_{0}$ the $(p,n$)-matrix whose $p\times p$ blocks are all null except for the last, which is the identity of ${\mathbb K}^p$ and denoted by $N_0$ the $(n\times n)$-matrix whose $p\times p$ blocks are all null except the sub-diagonals which are $Id_{{\mathbb K}^p}$. The previous formulas allow us to verify that the orbit of $(H_{0},N_{0})$ is open since the isotropy subalgebra is trivial.
\vskip 1mm
{\bf 2. Symplectic Reduction - Left invariant Symplectic foliation on $G_{n,p}({\mathbb K})$.} \\ 
Denote by $\omega^{+}=d\alpha^+$, $\alpha = (H,N)\in {\mathcal G}_{n,p}^{*}$ a left invariant symplectic form on $G_{n,p}$.\\
The action of $G_{n,p}({\mathbb K})$ on $G_{n,p}({\mathbb K})$ given by 
$$L_G : G_{n,p}({\mathbb K})\times  G_{n,p}({\mathbb K}) \longrightarrow G_{n,p}({\mathbb K}) , \qquad (\tau, \sigma) \longmapsto \tau\sigma \,( \mbox{ product in $G_{n,p}$ })$$
is a symplectic action. Morever $L_G$ is a Hamiltonian action;
a momentum mapping for the action is given by $$\mu: G_{n,p}({\mathbb K}) \longrightarrow  {\mathcal G}_{n,p}({\mathbb K})^{*} $$
$$ \sigma \longmapsto \mu (\sigma): x \longmapsto \langle \alpha^{+}(\sigma),x^{-}(\sigma) \rangle$$
where $ x \in {\mathcal G}_{n,p}$ and $x^{-}$ is the right invariant vector field associated to $x$.\\
The subgroup ${\mathcal H}:=M_{n,p}$ is (totally) isotropic for $\omega^{+}$ because  ${\mathcal H}$ is an abelian Lie group. Morever $L_{\mathcal H}: {\mathcal H} \times G_{n,p} \longrightarrow G_{n,p}$ is a hamiltonian action and a momentum mapping for $L_{\mathcal H}$ is given by:$$ m: G_{n,p} \longrightarrow  L({\mathcal H})^* $$ 
$$ \sigma \longmapsto \mu(\sigma)_{\arrowvert L({\mathcal H})}$$
where $L_{\mathcal H}$ is the Lie algebra of $\mathcal H$. We arrive at the following theorem :
\vskip 1mm
{\bf 2.1 Theorem .} Let $(G_{n,p},d((H,N)^+)$ defined as above. Then\\
{\bf 1.} $m^{-1}(H)$ is a closed subgroup of  $G_{n,p}$ and  ${\mathcal H} \subset m^{-1}(H)$\\
{\bf 2.} The canonical exact sequence of Lie groups 
$$\{\epsilon\} \rightarrow {\mathcal H} \rightarrow m^{-1}(H) \rightarrow  m^{-1}(H)/{\mathcal H} \rightarrow \{ \epsilon \} \leqno (2) $$ is split. It is also an exact sequence of affine Lie groups.\\
{\bf 3.} The reduced symplectic Lie group  $m^{-1}(H)/{\mathcal H}$ is isomorphic to :
 $$\begin{array}{cc} 
 G_{n-p,p}({\mathbb K}) & si\quad  \,\, n > p\\
 G_{p-1,1}({\mathbb K}) & si \quad n = p > 1\\
\end{array}$$
{\bf 4.} In the principal bundle
$$ m^{-1}(H) \stackrel{i}{\longrightarrow}  G_{n,p}({\mathbb K}) \stackrel{m}{\longrightarrow} \,  \Theta \leqno (3) $$
where $ \Theta$ is the set of matrices of rank $p$ in $M_{n,p}^{*}\equiv M_{p,n}$, the fiber is an affine Lie subgroup and $m$ is affine relative to the usual affine structure of $ \Theta \subset {\mathbb K}^{np}$.\\
We need the following lemma for which the proof is obvious.
\vskip 1mm
{\bf 2.2 Lemma .} The mapping $m$ is a surjective submersion.\\
Furthermore, 
$$ m((X,T))=H.T^{-1}  \, \mbox {\, where}  \quad (X,T)\quad \mbox {is \quad in }\quad G_{n,p}({\mathbb K})\leqno (4)$$ 
{\bf Proof.} The formula ({\bf 4}) is a direct consequence of the definition of $m$ and the Lemma 1.2. Thus it is clear that $m$ is a surjective submersion on the set  $\Theta$ of matrices of rank $p$ of $M_{n,p}(\mathbb K)$. 
\vskip 2mm
{\bf Proof of the theorem 2.1.}\\
Formula ({\bf 4}) implies that 
 $ m^{-1}(H) = \{(X,T) \in G_{n,p}({\mathbb K})/HT^{-1} = H \} $ is a (closed) subgroup of $G_{n,p}({\mathbb K})$ which contains ${\mathcal H}$. Morever the factor group $m^{-1}(H)/{\mathcal H}$ can be identified with $\{(0,T) \in G_{n,p}/HT^{-1} = H \}$. Thus ({\bf 2}) is an split sequence of Lie groups. 
On the other hand, since ${\mathcal H}$ is commutative and $\omega^+$ is exact, it follows that $L({\mathcal H}) \subset L({\mathcal H})^{\bot}$ and a straightforward shows that $L(m^{-1}(H))=L({\mathcal H})^{\bot}$.\\
Morever from the formula 
$$\omega (xy,z) = - \omega (y, [x,z]) $$ 
defining the left symmetric product in ${\mathcal G}_{n,p}({\mathbb K})$ (see formula ({\bf 1})) it turns out that $ L({\mathcal H})^{\bot}$ is a left symmetric subalgebra of ${\mathcal G}_{n,p}({\mathbb K})$ and that $L({\mathcal H})$ is a two-side ideal of $ L({\mathcal H})^{\bot}$. Consequently ({\bf 2}) is an exact sequence of affine Lie groups i.e $\mathcal H$, $m^{-1}(H)$, $m^{-1}(H)/{\mathcal H}$ are affine Lie groups and the applications of formula ({\bf 2}) are affine.\\
On the other hand a list of the elements of the matrix group $\{(0,T) \in G_{n,p}({\mathbb K})/ HT^{-1} = H \}$ allows to observe that the group $m^{-1}(H)/{\mathcal H}$ is isomorphic to the group $G_{n-p,p}({\mathbb K})$ if $n>p$ and isomorphic to $G_{n-1,1}({\mathbb K})$ if $n=p$. This proves statement {\bf 3.}\\
Now , we need to show that the manifold $\Theta$ is endowed with an affine structure which makes the momentum mapping affine.\\
Let $F$ be the subbundle of $TG_{n,p}({\mathbb K})$ tangent to the $L_{\mathcal H}$-orbit and $ F^{\bot}$ its symplectic orthogonal. Denote by $\mathcal F$ and ${\mathcal F}^{\bot}$ the associated foliations respectively. Since $\mathcal H$ is normal in $G_{n,p}({\mathbb K})$, the foliation ${\mathcal F}^{\bot}$ is defined either by the left invariant form on $G_{n,p}({\mathbb K})$ given by $\eta^{'}_{j}:=i(e^{+}_{j})\omega^+$ where $(e_{i})$'s form a basis of $L(\mathcal H)$ and $i$ denote the interior product, or by the closed forms (thus exact)  $\eta_{j}:=i(e^{-}_{j})\omega^+$. Obviously the $\eta_{j}$ are basic for the fibration. the forms $\bar{\eta}_{j}$ which are the projections of $\eta_{j}$ by $m$ define a local parallelism on $\Theta$. This parallelism is global and commutative because the $\bar{\eta}_{j}$ are exact. Hence $m$ is affine.  $\blacksquare$\\
Let $\alpha_1$ and $\alpha_2$ be the linear forms on  ${\mathcal G}_{n,p}({\mathbb K})$ defined respectively by $\alpha_{1}(x, u) = tr(H.x)$ and $\alpha_{2}(x, u) = tr(N.u)$. We get:\\
{\bf 2.3. Theorem.} $ker(d\alpha_{1})$ and $ker(d\alpha_{2})$ are supplementary symplectic Lie subalgebras of  $( {\mathcal G}_{n,p}({\mathbb K}) , d\alpha )$. So they determine two tranverse symplectic foliations left invariant on $ G_{n,p}({\mathbb K})$ \\
{\bf Proof.} Obviously $ker(d\alpha_{i})$ is a Lie subalgebra of ${\mathcal G}_{n,p}({\mathbb K})$. In addition the subspaces $ker(d\alpha_{i})$, $i = 1,2$ are in direct sum because ker$(d\alpha_{1})$$\bigcap$ker$(d\alpha_{2})$$=\{0\}$.\\
If we take $\alpha = (H_0,N_0)$ we directly observe that we have dim(ker$(\delta\alpha_{1}))$$=p^{2}(k-1)k$ and dim(ker$(\delta\alpha_{2}))$$=2p^{2}k$. We extends these results about the dimension for all $\alpha \in {\mathcal G}^{*}_{n,p}({\mathbb K})$ with open coadjoint orbit. consequently $ker(\delta\alpha_{1})$ and $ker(\delta\alpha_{2})$ are symplectic Lie subalgebras of $( {\mathcal G}_{n,p}({\mathbb K}) , d\alpha )$. $\blacksquare$
\vskip 2mm
Let $C(N_0)$ be the subalgebra of $gl({\mathbb K}^n)$ given by :
$$\left(\begin{array}{cccc}
A_0 &  &  &   \\
A_1 & \ddots &    &  \\
\vdots  & \ddots & \ddots   & \\
A_{k-1}  & \cdots & A_1 & A_0 \\
\end{array}\right) $$
{\bf 2.4. Scholie.} The Lie algebra $ker(d\alpha_{1})$ is isomorphic to the Lie algebra ${\mathcal G}_{n-p,p}({\mathbb K})$ if $n >p$ and ${\mathcal G}_{p-1,1}({\mathbb K})$ if $n=p$ while $ker(d\alpha_{2})$ is isomorphic to the semi-direct $M_{n,p}({\mathbb K})\rtimes C(N_0)$ if $n >p$ and to the semi-direct product $M_{p,p}({\mathbb K})\rtimes C(N_0)$ if $n=p$.\\
The following result is the infinitesimal version of Theorem 2.1. It gives a more precise and complete statement of Theorem 2.1.
\vskip 1mm  
{\bf 2.5 Proposition.} With the notations of Therem 2.1., if $I=L({\mathcal H})$, the canonical sequence of vectoriel spaces 
$$0 \rightarrow I \rightarrow I^{\bot} \rightarrow I^{\bot}/I \rightarrow 0 $$ is a split exact sequence of Lie algebras. It is also an exact sequence of left symmetric algebras.\\  
Furthermore the Lie algebra ${\mathcal G}_{n,p}(\mathbb K)$ decomposes as a direct sum of Lie subalgebras $ I^{\bot}$ and $C(N_{0})$ .
\vskip 2mm
{\bf 2.6 Theorem.} The symplectic Lie group $(G_{n,p},d\alpha^+)$ is endowed with two transversal left invariant symplectic foliations whose leaves are affine submanifolds of $G_{n,p}(\mathbb K)$. \\
\vskip 1mm
The following assertion specifies the number of open orbits in  ${\mathcal G}^{*}_{n,p}({\mathbb K})$ as well as the left invariant symplectic structures on $G_{n,p}({\mathbb K})$. 
\vskip 2mm
{\bf 2.7 Theorem.} \\ 
{\bf a.} There exist two open orbits of the coadjoint representation of $G_{n,p}({\mathbb K})$ if ${\mathbb K} ={\mathbb R}$ and only one if ${\mathbb K} ={\mathbb C} $\\
{\bf b.} Up isomorphism there is only one left invariant symplectic structure on $G_{n,p}$ i.e if $\omega$ and $\omega^{'}$ are two left invariant symplectic forms on $G_{n,p}({\mathbb K})$, then there exists an automorphism $\varphi$ of Lie algebra ${\mathcal G}_{n,p}(\mathbb K)$ such that :
$$ \omega_{\varepsilon}( ., . ) = {\omega}^{'}_{\varepsilon}(\varphi .,\varphi . ) . $$
\vskip 1mm
The following lemmas set up the main steps of the demonstration of Theorem 2.7. This lemma allow to count the $Ad^{*}_{G_{n,p}}$-orbits
\vskip 2mm
{\bf 2.8 Lemma.} If $Orb_{(H,M)}$ is the coadjoint orbit of $(H,M)\in {\mathcal G}_{n,p}^{*} \equiv M_{p,n}(\mathbb K) \times gl(n)$ then $Orb_{(H,M)}$ has an element $(H^{'}_{0},M^{'}_{0})\in {\mathcal G}_{n,p}^{*}$ such that $H^{'}_{0}= (0, \cdots , A_p) \in M_{p,n}$.\\
 Morever, if $Orb_{(H,M)}$ is open, then $H^{'}_{0}$ can be taken as $H^{'}_{0}= (0, \cdots , I_p)=H_0 $.
\vskip 2mm
{\bf Proof.} It suffices to remark that there exists $U\in GL_{o}({\mathbb K}^n)$ such that $HU^{-1} = H_0$ , that is obvious if we look at $H$ as the matrix of a mapping from ${\mathbb K}^n$ into ${\mathbb K}^p$ and $U^{-1}$ as the matrix of a change of basis in ${\mathbb K}^n$. Then the last formula of lemma 1.2 show the first assertion.\\
Now suppose that $Orb_{(H,M)}=Orb_{(H^{'}_{0},M^{'}_{0})}$ is a open orbit. This implies that we have $$(x,u)\in {\mathcal G}_{n,p}; \, ad^{*}_{(x,u)}(\alpha )= 0 \Rightarrow (x,u) =0$$
In particular $$ ad^{*}_{(x,0)}(H^{'}_0, M^{'}) =0 \Rightarrow x=0 $$
In other words
$$ tr(xH_{0}u)=0, u \in gl(n)\Rightarrow x=0 \leqno (**)$$
with $H_{0}= (0, A_p)$\\
But a straight calculation shows that (**) implies $A_p$ is invertible.\\
Finally there exists $U\in GL_{o}({\mathbb K}^n)$ such that $H^{'}_{0}U^{-1} = H_0$, but this relation fix only the last diagonal block of $U$ to the value $A^{-1}$. Therefore if $k \geq 2$ we able to take an other diagonal block egals to det$A^{-1}.Id_{{\mathbb K}^p}$ and completing the diagonale by the $1$, to construct a such matrix $U$ in $SL({\mathbb K}^n)$.
\vskip 2mm
{\bf 2.9 Lemma.} An open orbit $Orb_{(H_{0},M)}$ contains only one element $(H_{0},M^{'})$, where the block decomposition of $M^{'}$ can be written : $$M^{'}= \left(\begin{array}{cc} 
 M_1 & 0\\
 H_1 & 0  
\end{array}\right)$$ with $(H_{1},M_{1})\in {\mathcal G}_{n-p,p}^{*}$.
\vskip 2mm
{\bf Proof.} Because 
$$Ad^{*}_{(X, Id)}(H_0, M)=(H_0,M+X.H) $$
it is clear that there is only one $X \in M_{p,n}$ such that $M^{'}=M+MH_0$.
\vskip 2mm
{\bf 2.10 Lemma.} The linear forms $(H_{0},M^{'})$ and $(H_{o},P^{'})$ on ${\mathcal G}_{n,p}$ with $M^{'}=(H_{1},M_{1})$ and $P^{'}=(K_{1},P_{1})$, defined as in Lemma 2.9, belong to the same $Ad^{*}_{G_{n,p}}$-orbit if and only if $(H_{1},M_{1})$ and $(K_{1},P_{1})$ are in the same $Ad^{*}_{G_{n-p,p}}$-orbit.
\vskip 1mm
{\bf Proof.}
$(H_{0},M^{'})$ and $(H_{o},P^{'})$ are in the same orbit if and only if $\exists U \in GL_{o}({\mathbb K}^n)$ such that $UM^{'}U^{-1}=P^{'}$ and $H_{0}U^{-1} = H_o$.\\
The second condition implies $$U= \left(\begin{array}{cc} 
 U^{'} & 0\\
 X_1 & Id  
\end{array}\right) , \qquad with \qquad U_1 \in GL_0({\mathbb K}^{n-p})$$ 
On other hand  $UM^{'}U^{-1} = P^{'}$ in $gl({\mathbb K}^n)$ is equivalent to $Ad^{*}_{(X_1,U_1)}(H_1,M_1)=(K_1,P_1)$ in $ {\mathcal G}_{n-p,p}({\mathbb K})$
 $\blacksquare$\\  
{\bf Proof of the theorem 2.7.}\\
Following the previous Lemmas we have \\
cardinal$\left\{ open \qquad Ad^{*}_{G_{n,p}}-orbit \right\}$ = cardinal$\left\{ open \qquad Ad^{*}_{G_{p,p}}-orbit \right\}$\\
On the other hand, using similar arguments as those of the lemmas we can show that \\
cardinal$\left\{ open \qquad Ad^{*}_{G_{p,p}}-orbit \right\}$ = cardinal$\left\{ open \qquad Ad^{*}_{G_{1,1}}-orbit \right\}$ = $$\left\{\begin{array}{ll}
2 & \mbox{if} \quad {\mathbb K} = {\mathbb R}\\
1 & \mbox{if} \quad {\mathbb K} = {\mathbb C}
\end{array}\right.$$
If ${\mathbb K} = {\mathbb C}$ there is only one symplectic structure on $G_{n,p}$.\\
In the case ${\mathbb K} = {\mathbb R}$, it follows from the lemmas that every $Ad^{*}_{G_{n,p}}$-open orbit contains an element $(H_{0},N_{0})$ where
$$ N_{0}=\left(\begin{array}{ccccc}
0 &  &  &  &  \\
A_p & \ddots &  &   &  \\
  & I_p & \ddots &  & \\
  &  & \ddots & \ddots & \\
  &  &  & I_p & 0 \\ 
\end{array}\right) \quad \mbox{ with $A_p$  invertible}$$
Nevertheless two matrices as $N_0$ are conjugate by an element $P$ of $GL(n)$.\\
However the mapping $$ {\mathcal G}_{n,p}^{*} \longrightarrow  {\mathcal G}_{n,p}^{*},\qquad (g,M) \longmapsto ({P^t}^{-1}g,P^{-1}MP) $$ 
is dual of the mapping $$ \theta :{\mathcal G}_{n,p} \longrightarrow  {\mathcal G}_{n,p} \qquad (x,N) \longmapsto  ({P}^{-1}x,P^{-1}NP)$$
and these later is an automorphism of the Lie algebra ${\mathcal G}_{n,p}$.\\
Consequently if $\omega^+_1$ , $\omega^+_2$ are two left invariant symplectic forms on $G_{n,p}$, we have: 
 $$ \omega^{+}_{1} = \theta^{*}( \omega^{+}_{2}).\qquad \blacksquare$$
\vskip 2mm
{\bf Remark 2.}
Notice that the only one left invariant affine structure on $G_{n,p}$ is given by $(H_{0},N_{0})$.\\
The following result is an important consequence of the previous study.
\vskip 2mm
{\bf 2.11 Proposition.} The identity component of $G_{n,p}$ is diffeomorphic to an open $Ad^{*}_{G_{n,p}}$-orbit. Consequently $G_{n,p}$ and $\left\{ \alpha \in {\mathcal G}_{n,p}^{*} ; \, \alpha \mbox{ is Poisson-regular} \right\}$ are diffeomorphic.
\vskip 2mm
{\bf Proof.}(By induction)
We must prove that the orbital mapping in $(H_{0},N_{0})$
$$ (G_{n,p})_{o} \longrightarrow {\mathcal G}_{n,p}^{*}, \quad (X,U) \longmapsto (H_{0}U^{-1},UN_{0}U^{-1} + XH_{0}U^{-1})$$
has a trivial isotropy.\\
The result is obvious for $G_{1,1}$ and this implies that it is also true for $G_{p,p}$ ( thanks to a sequence double extension ).\\
Consider the case $G_{n,p}$ with $n=kp$, $k\geq 2$. The equality
$$ (H_{0}U^{-1},UN_{0}U^{-1} + XH_{0}U^{-1}) = (H_{0},N_{0})$$
implies that $UN_{0}U^{-1}$ and $N_0$ have the same p block-type ( in particular their last columns are zero). Hence $X=0$.\\
On the other hand, $U$ induces an element of $G_{n-p,p}$ belonging to the $Ad^{*}_{G_{n-p,p}}$-isotropy subgroup in $N_0$. Then, if the result is true in $G_{p,p}$ it is also true in $G_{2p,p}$ and $G_{kp,p}$ \, $p\geq 2$.$\blacksquare$
\vskip 2mm
{\bf 3. Left invariant Lagrangian foliations and Hess connection .}\\

{\bf 3.1.} Let's specify a little bit the Lie algebra isomorphisms indicated by the previous proposition. The isomorphism between Rad$(d\alpha_{1})$ and ${\mathcal G}_{n-p,p}(\mathbb K)$ is determined by the choice of a supplementary subspace of ${\mathcal K}(H_{0})=\left\{ X \in  M_{n,p}(\mathbb K), H_{0}.X=0 \right\}$ in $M_{n,p}(\mathbb K)$.\\
Denote by $X_0$, the element of $M_{n,p}(\mathbb K)$ formed by one column of zero-blocks except the last block which is $id_{{\mathbb K}^p}$. Then the mapping $ Rad(d\alpha_{1}) = \left\{ (0,u) \in {\mathcal G}_{n,p}(\mathbb K) , H_{0}.u=0 \right\} \longrightarrow   {\mathcal G}_{n-p,p}(\mathbb K) $ such that $(0,u) \longmapsto (U.u_{0}, \pi_{0}(u))$ defines such a isomorphism, where $\pi_{0}(u)$ denotes the matrix of the linear map given by $u$ restricted to ${\mathcal K}(H_{0})$.  Furthermore the image of the reduced symplectic form is the $2$-coboundary associated to $(H_{1},N_{1}) \in {\mathcal G}^{*}_{n-p,p}(\mathbb K)$ where $H_{1}$ and $N_{1}$ are defined by $tr(N_{0}u) = tr(H_{1}.u.X_{0}) + tr(N_{1}\pi_{0}(u))$ . We remark that $(H_{1},N_{1})$ has the same block type as $(H_{0},N_{0})$, so we can repeat the process of symplectic reduction in the same conditions. One obtains a decomposition of the space ${\mathcal G}_{n,p}(\mathbb K)$ as a direct sum of Lie subalgebras:
$${\mathcal G}_{n,p}(\mathbb K) = {\mathcal K}_{p-1}\oplus \cdots \oplus {\mathcal K}_0 \oplus C(N_{0})\oplus \cdots \oplus C(N_{p-1}) $$
where ${\mathcal K}_{i} $ comes from the totally isotropic ideal and ${\mathcal G}_{n-ip,p}(\mathbb K) = \{ K_{i}\times {\mathcal G}_{n-(i+1)p,p}(\mathbb K) \} \oplus C(N_{i}) $  for  $0 \leq i \leq p-1$ (orthogonal direct sum).
Using the canonical embedding of ${\mathcal G}_{n,p}({\mathbb K})$ in $gl({\mathbb K}^{n+p})$, the subspace $L= {\mathcal K}_{p-1}\oplus \cdots \oplus {\mathcal K}_0$ is identified to the subalgebra of strictly upper triangular $(p+1)\times (p+1)$-matrices by $p \times p$-blocks and $L^{'}=C(N_{0})\oplus \cdots \oplus C(N_{p-1})$ to the subalgebra of lower triangular $(p\times p$)-matrices by $(p\times p$)-blocks. As the  ${\mathcal K}_{i} $ and $C(N_{i})$ appear as the totally isotropic subspaces at every step of the successive reductions the subalgebras $L$ and $L'$ are lagrangian relative to $d\alpha$ where $\alpha \equiv (H_0,N_0)$ .\\
Denote by $\Lambda $ and ${\Lambda}^{'}$ the connected Lie subgroups of $G_{n,p}({\mathbb K})$ with Lie algebra $L$ and $L'$ respectively. The natural left actions of $\Lambda $ and ${\Lambda}^{'}$ on $(G_{n,p}({\mathbb K}),d\alpha^{+})$ being hamiltonnians, Theorem 3.1. of ([{\bf 3}]) allows us to assert that $\Lambda $ and ${\Lambda}^{'}$ are closed. So we have proved the following result:\vskip 2mm
{\bf 3.2 Theorem.} The symplectic Lie group $(G_{n,p},d\alpha^+)$ is endowed with two transversal left invariant lagrangian foliations with closed and affines leaves .
\vskip 2mm
Let $(G_{n,p}({\mathbb K}),\omega^{+})$ be endowed with its affine structure $\nabla$ defined by ({\bf 1}) \\
We recall that a connection on $(G_{n,p}({\mathbb K}),\omega^{+})$ is said to be {\sl symplectic} if and only if  $\nabla\omega= 0$ where $\omega := \omega^{+}_{\varepsilon}$, in other words:
$$ \nabla_{a}(\omega (b,c)) = \omega (\nabla_{a}b , c) + \omega (b, \nabla_{a}c) \quad \forall a=(x,u), b=(y,v), c=(z,r) \in {\mathcal G}_{n,p}({\mathbb K})$$
In that follows one identifies an element $a$ of ${\mathcal G}_{n,p}({\mathbb K})$ with $(x,u)$, where $x\in M_{n,p}$ and $u\in gl(n)$.\\
Let ${\mathcal L}$ and ${\mathcal L}^{'}$ be two lagrangian subalgebras of ${\mathcal G}_{n,p}({\mathbb K})$ such that ${\mathcal G} = {\mathcal L} \oplus {\mathcal L}^{'}$. Then we can write  $a = a_{1} + a_{2}$ where $a_{1}= (x_{1},u_{1}) \in {\mathcal L}$ and $a_{2}= (x_{2},u_{2}) \in {\mathcal L}^{'}$ . Then the left symmetric product on ${\mathcal G}_{n,p}({\mathbb K})$ is given by :
$$(x,0).(y,0)=(l,f)  \leqno (i)$$ where $l \in M_{n,p}$ is formed from a column of blocks such that the last belongs to $sl(p)$ and $f \in gl(n)$ is an element whose last line of blocks is zero.
$$(0,u).(0,v)= (l,-vu) \leqno (ii)$$ where $l$ is defined as above.
$$(x,0).(0,v)= (l_1,f_1)\leqno (iii)$$
$$(0,u).(y,0)= (l_2,f_2)\leqno (iv)$$
 where $l_1$,$l_2 \in M_{n,p}$, $f_1$,$f_2 \in gl(n)$ are defined as in (i).\\
The following result is a consequence of the previous discussion
\vskip 1mm
{\bf 3.3. Corollary.}
\textsl{There exists only one (torsion free ) left invariant symplectic connection $\tilde{\nabla}$  ( called the Hess'connection  ) such that }
$$ \tilde{\nabla}_{a}{\mathcal L} \subset {\mathcal L} , \,\, \tilde{\nabla}_{a'}{\mathcal L}^{'}\subset {\mathcal L}^{'}.$$ where $a \in {\mathcal L}$ et $a' \in {\mathcal L}^{'}$\\
{\sl This connection is defined by the products :}
$$ \tilde{\nabla}_{(x,0)^{+}}(y,0)^{+} = (l,0)^+ \quad \, ,  \quad  \tilde{\nabla}_{(0,u)^{+}}(0,v)^{+} = (0,-vu)^{+}$$
$$ \tilde{\nabla}_{(x,0)^{+}}(0,v)^{+} =0  \quad \, ,  \quad  \tilde{\nabla}_{(0,u)^{+}}(y,0)^{+} = (uy,0)^+$$
\vskip 1mm
Given the kind of results that we have described above, it is natural to ask ourself in what sense the group $G_{n,p}({\mathbb K})$ is not the symplectic double extension describes in ([{\bf 6}]).\\
The answer is clearly no if $k>1$; in fact to have $G_{n,p}({\mathbb K})$ be the symplectic double extension of $G_{n-p,p}({\mathbb K})$ in the sense of [{\bf 6}] it is necessary that $I^{\bot}$ be a Lie ideal of ${\mathcal G}_{n,p}({\mathbb K})$. Considering the proposition 2.3 involve the existence of a Lie ideal of $gl({\mathbb K}^n)$ isomorphic to ${\mathcal G}_{n-p,p}({\mathbb K})$.
\vskip 3mm
{\bf 4. $G_{n,p}(\mathbb K)$ \bf as  symplectic double extension of $G_{n-p,p}(\mathbb K)$} \\
We have shown in the previous paragraph that the techniquess of symplectic double extension developed in ([{\bf 5}]) do not apply to $G_{n,p}(\mathbb K)$ for $k > 1$.\\We reconsider the study of the canonical fibrations ({\bf 2}) and ({\bf 3} ) to try to understand how work this example. We have observed in 2.5. that for the symplectic Lie algebra  $({\mathcal G}_{n,p},d\alpha)$ where $\alpha \cong (H_{0},N_{0})$, a section of the canonical exact sequence
$$ 0 \rightarrow I \rightarrow I^{\bot} \rightarrow I^{\bot}/I \rightarrow 0 $$ is determined by the choice of an element $X_0$ in $M_{n,p}(\mathbb K)\backslash  {\mathcal K}(H_{0})$ where ${\mathcal K}(H_{0})$ is the kernel of $H_0$.\\
Conversely, we consider the reduced algebra  $({\mathcal G}_{n-p,p},d\alpha^{'})$ with $\alpha^{'} \equiv (H_{1},N_{1})$ where $H_{1}=(0,\cdots,0,I_{p})$ and 
$$ N_{1}=\left(\begin{array}{ccccc}
0 &  &  &  &  \\
I_p & \ddots &  &   &  \\
  & I_p & \ddots &  & \\
  &  & \ddots & \ddots & \\
  &  &  & I_p & 0 \\ 
\end{array}\right) $$
Let $i$ be the canonical inclusion from  $M_{n-p,p}$ to $M_{n,p}$ obtained by putting zero in the $n-p+1, \cdots , n$ rows, and let $Z \in M_{n,p}\backslash M_{n-p,p}$ have rank $p$ ( e.g.  $Z=H_{0})$).\\
Denote by $r: M_{n-p,p} \longrightarrow M_{n,p}$ the linear mapping given by 
$$r\circ i = id_{M_{n-p,p}} \mbox{ and } r(Z)=0 $$
and let $H \in M_{n,p}^{*}$ satisfy
$$ H\circ i =0 \mbox{ and } H.Z=Id_{p} \leqno (5)$$
We consider the regular representation  :
$$ \eta: {\mathcal G}_{n-p,p} \rightarrow gl(n),\, (x,u) \mapsto (i\circ u\circ r + H.i(x)) \leqno (6)$$
From $\eta$ we deduce a representation of Lie groups: \\
$$ \rho: G_{n-p,p} \rightarrow GL(n),\, (X,U) \mapsto (i\circ U\circ r + i(X).H) \leqno (7)$$
 verifying,
$$ \rho_{*,\epsilon} = \eta $$
Morever consider the inclusion :
$$R:M_{n,p}^{*} \times gl(n-p) \rightarrow gl(n) $$
given by 
$$ R(H,N)=i\circ N\circ r + Z.(H\circ r) \leqno (8) $$
Then we can state the following result.
\vskip 2mm
{\bf 4.1 Theorem.} Consider the symplectic Lie group $(G_{n-p,p},d(H_{1},N_{1})^{+})$ . If $H$ is the linear form given by ({\bf 4}), and $N=R(H_{1},N_{1})$ with $R$ given by ({\bf 7}), then $(G_{n,p},d(H,N)^{+})$ is a symplectic Lie group such that $(G_{n-p,p},d(H_{1},N_{1})^{+})$ is the reduced symplectic Lie group as described in theorem 2.1.
\vskip 2mm
{\bf Proof.} By definition of $H$ we have 
$$ m^{-1}(H)=M_{n,p} \times \rho (G_{n-p,p})$$
It remains to prouve that $d(H_{1},N_{1})^{+})$ is the reduced symplectic form of $d(H,N)^{+})_{\arrowvert m^{-1}(H)}$. Consider the decomposition $M_{n,p}({\mathbb K}) =i(M_{n-p,p}({\mathbb K}) )\otimes {\mathbb K}.Z$\\
On the other hand a straight verification proves that \\
$tr(H_{1}x)+tr(N_{1}u) = tr(R(H_{1},N_{1})\cdot (\eta (x,u))$  for all $(x,u)$ of ${\mathcal G}_{n-p,p}({\mathbb K})$.
Finally as $N=R(H_{1},N_{1})$ and $d((H,0)^+)$ vanish on $m^{-1}(H)$ it follows that $d((H_{1},N_{1}))^+)$ is the reduced form .This ends the proof. \qquad $\blacksquare$
\vskip 2mm
{\bf Remark 4.} A similar argument proves that for $p\geq 2$ : $G_{p,1}=Aff({\mathbb K}^{p})$ is a symplectic reduction of $G_{p,p}$ and $G_{p,p}$ is a symplectic double extension of $G_{p,1}$
\vskip 2mm
{\bf 4.2 Study of the fibration :}
$$ m^{-1}(H) \stackrel{i}{\longrightarrow} GA_{n,p}({\mathbb K}) \stackrel{m}{\longrightarrow} \, {\Theta } \leqno (3)$$ 
According to Theorem 1.2. we can consider the case where $\alpha \equiv (H_{0},N_{0})$. Let $ C=C(N_0)\bigcap GL({\mathbb K}^n)$ or, if we prefer $C$ is the Lie subgroup of $GL({\mathbb K}^n)$ formed by matrices of type 
$$\left(\begin{array}{cccc}
A_0 &  &  &   \\
A_1 & \ddots &    &  \\
\vdots  & \ddots & \ddots   & \\
A_{k-1}  & \cdots & A_1 & A_0 \\
\end{array}\right) $$  
with  $A_0$ is invertible.\\
According to Lemma 2.2. $m(C)=\{ (A_{k-1},\cdots,A_{1},A_{0}):  A_0 \mbox{ invertible }\}$, denote by $V_0$ this open set of $\Theta$.\\
Let be $V_{\gamma}=V_{0}.\sigma_{\gamma}$ where $\gamma =(i_{1},\cdots,i_{p})$ and $\sigma_{\gamma}$ is the element of $GL({\mathbb K}^n)$ which realizes the permutation of the $k$ last columns with the columns indexed by $\gamma$ in $M_{p,n}({\mathbb K})$. The sets $V_{\gamma}$ are clearly the open set of $\mathcal \Theta$ and $m$ defines a diffeomorphism of $C.\sigma_{\gamma}$ on $V_{\gamma}$. Denote $S_{\gamma}$ the embedding of $V_{\gamma}$ in $GL({\mathbb K}^n)$ such that $S_{\gamma}(V_{\gamma})=C.\sigma_{\gamma}$ and $m\circ S_{\gamma}=id_{U_{\gamma}}$ for all multi-indices $\gamma$.\\
We have the following result:
\vskip 2mm
{\bf 4.3 Lemma.} The $V_{\gamma}$ as defined above , with $\gamma =(i_{1},\cdots,i_{p})$ for $1\leq i_{1}<i_{2}<\cdots <i_{p}\leq n$, form an open trivializing cover for the fibration ({\bf 2}).
\vskip 1mm
{\bf Proof.}\\
The elements $V_{\gamma}$ recover $\Theta$ since $\Theta$ is formed by matrices of rank $p$ in $M_{p,n}({\mathbb K})$ whose first $p$ columns are independent in ${\mathbb K}^{p}$. By contrast $V_{0}$ is formed of matrices for which the last $p$ colums are independent. \\
To show that the $V_{\gamma}$ are trivialising for the fibration ({\bf 3}) amounts to proving that we have $m^{-1}(V_{\gamma})=S_{\gamma}(V_{\gamma}).m^{-1}(H_{0})$. For all $\sigma$ of $G_{n,p}({\mathbb K})$ we have the identity $m^{-1}(m(\sigma))=\sigma .m^{-1}(H_{0})$. Indeed if $\sigma^{'} \in m^{-1}(m(\sigma))$ we have by definition of $m$ the formula,
$$\langle Ad_{\sigma}^{*}(H,0);(X,0) \rangle = \langle Ad_{\sigma^{'}}^{*}(H,0);(X,0) \rangle \mbox{ for all } X \mbox{ in } M_{n,p}({\mathbb K})$$
or what amounts to the same thing
$$\langle Ad_{\sigma^{-1}\sigma^{'}}^{*}(H,0);(X,0) \rangle = \langle (H,0);(X,0) \rangle \mbox{ for all } X \mbox{ in } M_{n,p}({\mathbb K}) $$
This last relation means that $\sigma^{-1}\sigma^{'} \in m^{-1}(H_{0})$ by the identity below and the fact that the $V_{\gamma}$ are open and trivialising.The trivialisations are then given by the map 
$\phi_{\gamma}: m^{-1}(V_{\gamma}) \rightarrow V_{\gamma} \times m^{-1}(H_{0}); \, S_{\gamma}(\alpha).\sigma \mapsto (\alpha, \sigma)$ for all multi-indices $\gamma =(i_{1},\cdots,i_{p})$ with $\leq i_{1}<i_{2}<K <i_{p}\leq n$. The cocycle defining the fibration ({\bf 2}) is given by
$$ {\Gamma}_{\gamma_{1}\gamma_{2}} : V_{\gamma_{1}}\bigcap V_{\gamma_{2}} \rightarrow m^{-1}(H_{0}); \,\alpha \mapsto S_{\gamma_{1}}^{-1}(\alpha).S_{\gamma_{2}}(\alpha) \leqno (9)$$
Thus we have proved the following result.
\vskip 2mm
{\bf 4.4 Proposition.}
The manifold $G_{n,p}({\mathbb K})$ is diffeomorphic to $\coprod V_{\gamma} \times m^{-1}(H_{0})/\sim $ where $(\alpha,\sigma)\sim (\beta,\tau)$ if and only if $\alpha = \beta$ and $\sigma = \Gamma_{\gamma_{1}\gamma_{2}}(\tau)$ for $(\alpha,\sigma)$ in  $V_{\gamma_{1}} \times m^{-1}(H_{0})$ and $(\beta,\tau)$ in $V_{\gamma_{2}} \times m^{-1}(H_{0})$. The cocycle $\Gamma_{\gamma_{1}\gamma_{2}}$ is defined by ({\bf 9}).
\vskip 1cm
{\it Acknowledgements.} The third author was partially supported by NSF grant DMS-0204100 and The Department of Mathematics of Penn state University. He also wishes to thank The Department of Romances Languages and the Department of Mathematics at The University of North Carolina at Chapel Hill for their kind welcome during his stay.
Thanks to Jim Stasheff for editing this paper and thanks to  Patrick Eberlein for his helpful discussions and suggestions.
\vskip 1cm

$\begin{array}{ll}
Jean Michel \, DARDIE &  Alberto \, MEDINA \\
 Universit\mbox{\'e} \, de \, Montpellier 2 & Universit\mbox{\'e} \, de \, Montpellier 2\\ 
D\mbox{\'e}partement \, de \, Math\mbox{\'e}matiques & D\mbox{\'e}partement \, de \, Math\mbox{\'e}matiques  \\
34095 \, Montpellier \, cedex 5 & 34095 \, Montpellier \, cedex 5 \\
UMR \, 5149 \, du \, CNRS & UMR \, 5149 \, du \, CNRS \\
 dardie@darboux.math.univ-montp2.fr & medina@darboux.math.univ-montp2.fr 
\end{array}$
\vskip 5mm
$$\begin{array}{l}
 Hassene \, SIBY \\
 Universit\mbox{\'e} \, de \, Montpellier 2\\
 D\mbox{\'e}partement \, de \, Math\mbox{\'e}matiques\\
 34095 \, Montpellier cedex 5\\
 UMR \, 5149 \, du \, CNRS\\
 siby@darboux.math.univ-montp2.fr
\end{array}$$


\newpage
\begin{thebibliography}{99}

\bibitem[1]{1}\textsc{R. Abraham }  \textsc{J.E. Marsden,} Fundation of Mechanics, \textit{Benjamin/Cummings Pub., Princeton, N.J,}  1978.

\bibitem[2]{2} \textsc{Y. Agoaka ,} Uniqueness of Left Invariant Symplectic Structures on the Affine Lie Group, \textit{Proc.Amer.Math.Soc. } \textbf{Vol.129},$N^{\circ }$9 (2001), pp 2753-2762.


\bibitem[3]{3} \textsc{M. Bordemann ,}   \textsc{A. Medina, }  \textsc{A. Ouadfel,}  Le groupe affine comme vari\'et\'e symplectique, \textit{Tohuku Math. J. } \textbf{45},$N^{\circ }$3 (1993), pp 423-436.


\bibitem[4]{} \textsc{J.M. Dardi\'e ,} Groupes de Lie \`a structure symplectique ou k$\ddot{a}$hl\'erienne et double extension. \textit{ Th\`ese de Doctorat, Univ.Montpellier II },1993.


\bibitem[5]{5} \textsc{J.M. Dardi\'e,}   \textsc{A. Medina,} Groupes de Lie \`a structure symplectique invariante, \textit{in ``S\'eminaire Gaston Darboux'', Montpellier } \textbf{t.328, s\'erie I},(1990-1991),p.77-85.


\bibitem[6]{6}\textsc{J.M. Dardi\'e,}  \textsc{A. Medina,} Double extension symplectique d'un groupe de Lie symplectique, \textit{Adv. Math. } \textbf{117} (1996),
pp 208-227.


\bibitem[7]{7}\textsc{J.M. Dardi\'e,} \textsc{A. Medina,} Alg\`ebres de Lie k$\ddot{a}$hl\'eriennes et double extension , \textit{J. Algebra. } \textbf{185} (1996),
pp 774-795.


\bibitem[8]{8}\textsc{H.Hess}, Connxions on Symplectic Manifolds and Geometric Quantization , Differential geometrical methods in mathematical physics, \textit{( Proc. Conf., Aix-en-Provence/Salamanca )}  (1979),
pp 153-166.

\bibitem[9]{9}\textsc{A. Medina,} \textsc{Ph. revoy} Groupes de Lie \`a structure symplectique invariante , \textit{Groupoids and Integrable Systems, M.R.S.I.Pub. Berkeley}  (1991),
pp 247-266.



\bibitem[10]{10}\textsc{M. Rais, } La repr\'esentation coadjointe du groupe affine,  \textit{Ann. Inst. Fourier} \textbf{28} (1978),
pp 207-237.



\bibitem[11]{11}\textsc{A. Weinstein, } Symplectic Manifolds and Their Lagrangian submanifolds,  \textit{Adv. Math. } \textbf{6} (1971),
pp 329-346.

\end{thebibliography}
\end{document}